\documentclass[draft]{amsart}

\usepackage{enumerate}
\usepackage{amsmath,amssymb,epsfig,graphicx,url}
\usepackage{color}
\usepackage{todonotes}
\usepackage[all,knot,arc]{xy}

\newtheorem{thm}{Theorem}[section]
\newtheorem{lem}[thm]{Lemma}
\newtheorem{cor}[thm]{Corollary}

\newtheorem{exam}[thm]{Example}
\newtheorem{defn}[thm]{Definition}
\newtheorem{problem}{Problem}

\newenvironment{pf}{{\bf Proof}}

\renewcommand{\leq}{\leqslant}
\renewcommand{\geq}{\geqslant}

\renewcommand{\ker}{{\rm ker}}

\newcommand{\aut}{\mathop{\mathrm{Aut}}}
\newcommand{\Aut}{\aut}

\newcommand{\rank}{\operatorname{rank}}

\newcommand{\genset}[1]{\ensuremath{\langle\: #1 \:\rangle}}
\newcommand{\trans}{\mathcal{T}_{n}}
\newcommand{\sym}{\mathcal{S}_{n}}
\newcommand{\alt}{\mathcal{A}_{n}}

\newcommand{\agl}{\mbox{\rm AGL}}
\newcommand{\Sz}{\mbox{\rm Sz}}
\newcommand{\psl}{\mbox{\rm PSL}}

\newcommand{\pgl}{\mbox{\rm PGL}}

\newcommand{\pgaml}{\mbox{\rm P}\Gamma {\rm L}}
\newcommand{\agaml}{\mbox{\rm A}\Gamma {\rm L}}
\newcommand{\psigl}{\mbox{\rm P}\Sigma {\rm L}}
\newcommand{\asigl}{\mbox{\rm A}\Sigma {\rm L}}

\newcommand{\dihed}[1]{\ensuremath{D_{#1}}}

\renewcommand{\to}{\longrightarrow}

\begin{document}



\title[Groups and Transformation Semigroups]{The Classification of Partition Homogeneous Groups with Applications to Semigroup Theory}
\author{Jorge Andr\'{e}}
\author{Jo\~{a}o Ara\'{u}jo}
\author{Peter J. Cameron}

\address[Andr\'{e}]
{Departamento de Matem\'{a}tica,
Universidade Nova de Lisboa\\}
\email{\url{jmla@fct.unl.pt}}

\address[Ara\'{u}jo]
{Universidade Aberta
and
Centro de \'{A}lgebra \\
Universidade de Lisboa \\
Av. Gama Pinto, 2, 1649-003 Lisboa \\ Portugal}
\email{\url{jaraujo@ptmat.fc.ul.pt}}

\address[Cameron]{Mathematical Institute\\ University of St Andrews\\
St Andrews, Fife KY16 9SS\\Scotland}
\email{\url{pjc@mcs.st-andrews.ac.uk }}


\begin{abstract}
Let $\lambda=(\lambda_1,\lambda_2,\ldots)$ be a \emph{partition} of $n$, a
sequence of positive integers in non-increasing
  order with sum $n$. Let
$\Omega:=\{1,\ldots,n\}$. An ordered
partition $P=(A_1,A_2,\ldots)$ of $\Omega$ has \emph{type} $\lambda$ if
$|A_i|=\lambda_i$.

Following Martin and Sagan, we say that $G$ is \emph{$\lambda$-transitive} if,
for any two ordered partitions $P=(A_1,A_2,\ldots)$ and $Q=(B_1,B_2,\ldots)$ of
$\Omega$ of type $\lambda$, there exists $g\in G$ with $A_ig=B_i$ for all $i$. A group 
$G$ is said to be \emph{$\lambda$-homogeneous} if, given two ordered partitions $P$ and $Q$ as above, inducing the sets $P'=\{A_1,A_2,\ldots\}$ and $Q'=\{B_1,B_2,\ldots\}$, there
exists $g\in G$ such that $P'g=Q'$. Clearly a $\lambda$-transitive group is
$\lambda$-homogeneous.

The first goal of this paper is to classify the $\lambda$-homogeneous groups (Theorems \ref{lambdahom} and \ref{lambdatrans}). The second goal is to apply this classification to a problem in semigroup theory.

Let $\trans$ and $\sym$ denote the transformation monoid and the symmetric group on $\Omega$, respectively. Fix a group $H\leq \sym$. Given a non-invertible transformation $a\in \trans\setminus \sym$  and a  group $G\leq \sym$, we say that $(a,G)$ is an \emph{$H$-pair}  if the semigroups generated by $\{a\}\cup H$ and $\{a\}\cup G$ contain the same non-units, that is, $\langle a,G\rangle\setminus G=\langle a,H\rangle\setminus H$.
Using the classification of  the $\lambda$-homogeneous groups we classify all the $\sym$-pairs (Theorem \ref{sympairs}).  For a multitude of transformation semigroups this theorem immediately implies a description of their automorphisms, congruences, generators and other relevant properties  (Theorem \ref{illustr}). 

This topic involves both group theory and semigroup theory; we have attempted to include enough exposition to make the paper self-contained for researchers in both areas.

 The paper finishes with a number of open problems on permutation and linear groups. 
\end{abstract}
\maketitle
\medskip

\noindent{\em Date:} 5 April 2013\\
{\em Key words and phrases:} Transformation semigroups, permutation groups, primitive
groups, $\lambda$-transitive, GAP\\
{\it 2010 Mathematics Subject Classification:}  20B30, 20B35,
20B15, 20B40, 20M20, 20M17. \\
{\em Corresponding author: Jo\~{a}o Ara\'{u}jo}





\newpage

\section{Introduction}
For notation and basic results on group theory we refer the reader to \cite{cam,dixon}; for semigroup theory we refer the reader to \cite{Ho95}.
A permutation group $G$ on $\Omega:=\{1,\ldots,n\}$ (a subgroup of $\sym$) is said
to be $t$-transitive if it acts transitively on the set of $t$-tuples of
distinct elements of $\Omega$, and is $t$-homogeneous if it acts
transitively on the set of $t$-element subsets of $\{1,\ldots,n\}$. Clearly
a $t$-transitive group is $t$-homogeneous. Livingstone and Wagner~\cite{lw}
showed that the converse is true if $5\le t\le n/2$; the $t$-homogeneous, but not $t$-transitive
groups for $t=2,3,4$ were determined by
Kantor~\cite{kantor,kantor2}.

Martin and Sagan~\cite{ms} defined the notion of $\lambda$-transitivity of a
permutation group, as a generalization of $t$-transitivity, where $\lambda$
is a partition of $n$ (the degree of the permutation group). In the case where
$\lambda=(n-t,1,1,\ldots,1)$, $\lambda$-transitivity is just $t$-transitivity.
There is a  weakening of this notion, here called
$\lambda$-homogeneity; this notion is natural in permutation groups and also in transformation
semigroups. Our first goal is to investigate the relationship between $\lambda$-transitivity and $\lambda$-homogeneity.

Let $\lambda=(\lambda_1,\lambda_2,\ldots)$ be a \emph{partition} of $n$, a
sequence of positive integers in non-increasing order with sum $n$. Now let
$G$ be a  permutation group on $\Omega$, where $|\Omega|=n$. An ordered
partition $P=(A_1,A_2,\ldots)$ of $\Omega$ has \emph{type} $\lambda$ if
$|A_i|=\lambda_i$.

Following Martin and Sagan, we say that $G$ is \emph{$\lambda$-transitive} if,
for any two partitions $P=(A_1,A_2,\ldots)$ and $Q=(B_1,B_2,\ldots)$ of
$\Omega$ of type $\lambda$, there exists $g\in G$ with $A_ig=B_i$ for all $i$. A group 
$G$ is said to be \emph{$\lambda$-homogeneous} if, given two ordered partitions $P$ and $Q$ as above, inducing the sets $P'=\{A_1,A_2,\ldots\}$ and $Q'=\{B_1,B_2,\ldots\}$, there
exists $g\in G$ such that $P'g=Q'$. Clearly a $\lambda$-transitive group is
$\lambda$-homogeneous.

Our motivating question is whether the analogue of the Livingstone--Wagner
theorem holds for partition homogeneity and transitivity. Indeed, we see that
a group of degree $n$ which is $t$-homogeneous but not $t$-transitive is 
$\lambda$-homogeneous but not $\lambda$-transitive, where $\lambda$ is the
partition $(n-t,1,1,\ldots,1)$ (with $t$ ones).

Our first theorem characterizes the  $\lambda$-homogeneous groups.
The partition 
$\lambda=(1,1,\ldots,1)$ is excluded since every permutation group is
$\lambda$-homogeneous but only the symmetric group is $\lambda$-transitive.
 
\begin{thm}\label{lambdahom}
Let $\lambda$ be a partition of $n$ (other than $(1,1,,\ldots,1)$). Let
$G$ be a $\lambda$-homogeneous subgroup of $\sym$. If $G$ is not 
$\lambda$-transitive, then one of the following holds:
 \begin{enumerate}
 \item $G$ is intransitive, and either
  \begin{enumerate}
  \item $\lambda=(k,k,\ldots,k)$ ($t$ parts, with $tk=n$),  $G$ fixes a point
  and acts as $\mathcal{S}_{n-1}$ or $\mathcal{A}_{n-1}$ on the remaining
  points; or
  \item $\lambda=(k,k)$ with $k=3$ or $5$, and $G$ fixes a point and acts
  as $\agl(1,5)$ ($n=6$) or $\pgl(2,8)$ or $\pgaml(2,8)$ ($n=10$)
  on the remaining $5$ or $9$ points.
  \end{enumerate}
 \item\label{trans} $G$ is transitive and one of the following occurs:
  \begin{enumerate}
  \item $\lambda=(n-t,1,1,\ldots,1)$, and $G$ is $t$-homogeneous but not
  $t$-transitive;
  \item $n=6$, $\lambda=(3,3)$, $G=\psl(2,5)$;
  \item $\lambda=(3,2,1,\ldots,1)$, $G$ is the Mathiu group $M_n$ ($n=12$ or
  $n=24$);
  \item $\lambda=(2,2,1,\ldots,1)$, and either $G$ is the Mathieu group $M_n$
  ($n=11,12,23,24$) or $G=\pgaml(2,8)$ ($n=9$).
  \end{enumerate}
 \end{enumerate}
\end{thm}

The exceptional groups listed in (1)(b) and (2)(b) of the previous theorem all
have two orbits on $k$-subsets of $\{1,\ldots,n\}$, when $(k,n)$ is either
$(3,6)$ or $(5,10)$; the sets in one orbit are the complements of those in
the other.

The previous result classifies $\lambda$-homogeneous groups modulo a classification of $\lambda$-transitive groups. 
Our next result is the classification of $\lambda$-transitive groups.
In order to state the theorem, we need some preliminaries.

Clearly $\sym$ is $\lambda$-transitive for every partition $\lambda$ of $n$,
and $\alt$ is $\lambda$-transitive for every partition except
$\lambda=(1,1,\ldots,)$. So we may exclude these two groups.

If $G$ is $\lambda$-transitive, and $t$ is any integer with $0<t<n$ which is
a sum of some of the parts of $\lambda$, then $G$ is $t$-homogeneous. In
particular, if $G$ is $\lambda$-transitive with $\lambda\ne(n)$, then $G$
is transitive.

We say that the pair $(G,\lambda)$ is \emph{standard} if
\begin{itemize}
\item $G$ is $t$-homogeneous, with $t<n/2$;
\item $\lambda=(n-t,k_1,\ldots,k_{r-1})$;
\item the stabiliser of a $t$-set in $G$ acts
$(k_1,\ldots,k_{r-1})$-transitively on it.
\end{itemize}

It is clear that, if $(G,\lambda)$ is standard, then $G$ is
$\lambda$-transitive. Conversely, if the largest part of $\lambda$ is
$n-t$, with $t<n/2$, and $G$ is $\lambda$-transitive, then $G$ is standard.
(For with this assumption, $\lambda$ has only one part of size $n-t$; so $G$
must act $(n-t)$-homogeneously, hence $t$-homogeneously, and the third
hypothesis then implies $\lambda$-transitivity.)
Later (in Subsection \ref{eti}) we will prove that if $G$ is $\lambda$-transitive then the largest part of $\lambda$ is
$n-t$, with $t\le n/2$. Therefore we have the following result.
\begin{thm}\label{lambdatrans}
Let $\lambda\neq(n)$ be a partition of $n$, and $G\le\sym$, with $G\ne\alt,\sym$.
Then  $G$ is $\lambda$-transitive if and only if $(G,\lambda)$ is standard.
\end{thm}

Therefore, to check whether a given group G is $\lambda$-transitive, we  only have to check
that $G$ is $t$-homogeneous (where the largest part of $\lambda$ is $n-t$) and
that the stabiliser of a $t$-set is $\lambda'$-transitive, where $\lambda'=(k_1,\ldots,k_{r-1})$, when $\lambda=(n-t,k_1,\ldots,k_{r-1})$. 

To inspect the groups and partitions regarding standardness we consider two cases: first the $t$-transitive groups and then  the $t$-homogeneous, but not $t$-transitive groups. 

If $G$ is
$t$-transitive, the stabiliser of a $t$-set induces $\mathcal{S}_{t}$ on it and so is
$\lambda'$-transitive for any partition $\lambda'$ of $t$.

If $G$ is
$t$-homogeneous but not $t$-transitive, we only have to look at Kantor's list (given later)
of such groups for $t=2, 3, 4$, and in each case figure out what is the group
induced on a $t$-set by its setwise stabiliser (here called $H$), and check for which partitions
$\lambda'$ this group acts $\lambda'$-transitively.  Observe that the group  $H$ is
well-defined since all $t$-sets are equivalent under the action of the group
$G$. 

Now, if $t=2$ then $G\le\asigl(1,q)$ and $H$ is the trivial group 
 (so that $\lambda'$ has to be $(2)$). For $t=3$, $H$ is
the cyclic group of order $3$ (so that $\lambda'$ can be either $(3)$ or $(2,1)$), except when $G$ is $\agl(2,8)$ or $\agaml(2,32)$ in which cases $H$ it is trivial
(so $\lambda'$ can only be $(3)$). Finally, if $t=4$, then $H$ is the alternating group on $4$ points 
for $\pgaml(2,8)$ (so $\lambda'$ is $(4)$, $(3,1)$, $(2,2)$ or $(2,1,1)$), and the
Klein group in the other two cases (so only $(4)$ and $(3,1)$ are possible).

%

To summarize, if $G$ is $\lambda$-transitive, $H$ denotes the group induced on a $t$-set by its stabiliser, and $V_{4}$ is the Klein group, we have the following:
\begin{itemize}
\item $t=2$, $G\le\asigl(1,q)$, $H=\{1\}$: $\lambda=(n-2,2)$;
\item $t=3$, $G\le\psigl(2,q)$, $H=\mathcal{A}_3$: $\lambda=(n-3,3)$ or
$(n-3,2,1)$;
\item $t=3$, $G=\agl(1,8)$ or $\agaml(1,32)$, $H=\{1\}$: $\lambda=(n-3,3)$;
\item $t=3$, $G=\agaml(2,8)$, $H=\mathcal{A}_3$: $\lambda=(n-3,3)$ or
$(n-3,2,1)$;
\item $t=4$, $G=\pgl(2,8)$ or $\pgaml(2,32)$, $H=V_4$: $\lambda=(n-4,4)$ or
$(n-4,3,1)$;
\item $t=4$, $G=\pgaml(2,8)$, $H=\mathcal{A}_4$: $\lambda$ is any partition
with largest part $n-4$ except $(n-4,1,1,1,1)$.
\end{itemize}

The classification of $\lambda$-homogeneous groups is a very natural problem in group theory, and it is quite a mystery why it took so long for this idea to appear in the literature, after the introduction of $\lambda$-transitive groups by Martin and Sagan in  \cite{ms}. However,  the purpose of the present paper is not just to generalize $\lambda$-transitivity, but to study a very natural  problem in semigroup theory that we now explain.

 Let $\trans$ denote the  monoid consisting of mappings from $\Omega$ to $\Omega$. The monoid $\trans$ is usually called the full transformation semigroup. In \cite{lm}, Levi and McFadden proved the following result.

\begin{thm}\label{original}
Let $a\in \trans\setminus \sym$. Then
\begin{enumerate}
\item  $\langle g^{-1}ag\mid g\in \sym\rangle$ is idempotent generated;
\item $\langle g^{-1}ag\mid g\in \sym\rangle$ is regular.
\end{enumerate}
\end{thm}

Using a beautiful argument, McAlister \cite{mcalister} proved that the semigroups $\langle g^{-1}ag\mid g\in \sym\rangle$ and $\langle a, \sym\rangle\setminus \sym$ (for $a\in \trans\setminus \sym$) have exactly the same set of  idempotents; therefore, as $\langle g^{-1}ag\mid g\in \sym\rangle$ is idempotent generated, it follows that
\begin{eqnarray}\label{(1)}
 \langle g^{-1}ag\mid g\in \sym\rangle=\langle a, \sym\rangle\setminus \sym.
\end{eqnarray}

In another direction, Levi \cite{levi96}  proved that, for every $a\in \trans\setminus \sym$, we have
\begin{eqnarray}\label{(2)}
\langle g^{-1}ag\mid g\in \sym\rangle\setminus \sym = \langle g^{-1}ag\mid g\in \alt\rangle\setminus \alt  .
\end{eqnarray}

%
%
%

The two next theorems appear in \cite{ArMiSc}  and generalize 
Theorem \ref{original}. 

\begin{thm}\label{t3}
If $n\geq 1$ and $G$ is a subgroup of $\sym$, then the following are equivalent:
\begin{enumerate}
\item The semigroup
$\genset{g^{{-1}}ag \mid g\in G}$ is idempotent generated
for all $a\in \trans\setminus\sym$.
\item One of the following is valid for $G$ and $n$:
\begin{enumerate}
\item $n=5$ and $G$ is $\agl(1,5)$;
\item $n=6$ and $G$ is  $\psl(2,5)$ or $\pgl(2,5)$;
\item $G$ is $\alt$ or $\sym$.
\end{enumerate}
\end{enumerate}
\end{thm}

\begin{thm}\label{th2}
If $n\geq 1$ and $G$ is a subgroup of $\sym$,  then the following are equivalent:
\begin{enumerate}
    \item The semigroup $\genset{g^{{-1}}ag \mid g\in G}$ is regular for all $a\in \trans\setminus\sym$.
\item One of the following is valid for $G$ and $n$:
    \begin{enumerate}
    \item $n=5$ and $G$ is  $ C_5,\ \dihed{5},$ or $\agl(1,5)$;
    \item $n=6$ and $G$ is $ \psl(2,5)$ or $\pgl(2,5)$;
    \item $n=7$ and $G$ is $\agl(1,7)$;
    \item $n=8$ and $G$ is $\pgl(2,7)$;
    \item $n=9$ and $G$ is  $\psl(2,8)$ or $\pgaml(2,8)$;
    \item $G$ is $\alt$ or $\sym$.
\end{enumerate}
\end{enumerate}
\end{thm}

%
%

The groups  $G\leq \sym$ possessing the property exhibited by $\sym$ in (\ref{(1)}) above are classified in \cite{acmn} as follows. 

\begin{thm}\label{mainold}
If $n\geq 1$ and $G$ is a subgroup of $\sym$,  then the following are equivalent:
\begin{enumerate}
    \item for all $a\in \trans \setminus \sym$ we have $$ \langle g^{{-1}}ag\mid g\in G\rangle = \langle a,G\rangle \setminus G  ;$$
\item one of the following is valid for $G$ and $n$:
    \begin{enumerate}
    \item $n=5$ and $G$ is $ \agl(1,5)$;
    \item $n=6$ and  $G$ is  $ \psl(2,5)$ or $\pgl(2,5)$;
    \item $n=9$ and  $G$ is  $\psl(2,8)$ or $\pgaml(2,8)$;
    \item  $G$ is $\{1\}$, $\alt$ or $\sym$.
\end{enumerate}
\end{enumerate}
\end{thm}

Within the boundaries of this study on how the properties of the group shape the structure of the semigroup, and also  along the lines for future investigations suggested by Levi, McAlister and McFadden in \cite{lmm},  the {\em ultimate goal}  is to classify all the pairs $(a,G)$ (where $a$ is a singular map and $G$ is a group of permutations) such that the semigroup $\langle a,G\rangle \setminus G$ has a given property $P$. As particular instances of this general problem we have the following list (in what follows, $a$ will be always a transformation in $\trans\setminus \sym$ and $G$ will be a subgroup of $\sym$):

\begin{enumerate}
\item Classify the pairs $(a,G)$ such that the semigroup $\langle a,G\rangle \setminus G$  is regular. 
\item Classify the the pairs $(a,G)$  such that the semigroup $\langle a,G\rangle \setminus G$  is generated by its idempotents.
\item Classify the the pairs $(a,G)$  such that   $\langle a,G\rangle\setminus G=\langle g^{-1}ag\mid g\in G\rangle$. 
\end{enumerate}
 
In \cite{ArCa12} some progress has been  made on this direction and the results therein will certainly be crucial ingredients in the solution of the  {\em ultimate goal}.

In this setting, suppose we have a group $H\leq \sym$. A pair $(a,G)$, where $G\leq \sym$ and $a\in \trans\setminus \sym$,  such that  
 $\langle a,G\rangle\setminus G=\langle a,H\rangle\setminus H$,  is called an $H$-pair. As the semigroups $\langle a,\sym\rangle\setminus\sym$ have been intensively studied under the name of \emph{$S_{n}$-normal semigroups} (see \cite{lm}, the references therein and also the papers that cite this one), it is natural to ask for a classification of the $\sym$-pairs. 
 
 By (\ref{(2)}) above we already know that $(a,\alt)$ is an $\sym$-pair, for every $a\in \trans\setminus \sym$; and in \cite{mcalister} it is proved that if $e^{2}=e\in \trans$ and $|\Omega e|=n-1$, then $(e,G)$ is an $\sym$-pair if and only if $G$ is $2$-homogeneous. 
 
 Our second main result dramatically  generalizes these two results. Using our results about $\lambda$-homogeneous groups, we  provide the full classification of the $\sym$-pairs (Theorem \ref{sympairs}). Therefore we provide one instance of the {\em ultimate goal}:  we classify all pairs $(a,G)$ that generate a semigroup with a given property $P$ (where the property $P$ is in this case  {\em `$(a,G)$  generates the semigroup $\langle a,\sym\rangle\setminus\sym$'}).  This classification immediately implies a large number of results about semigroups. 
In fact, since \emph{almost} everything is known about the semigroups  
$\langle g^{-1}ag\mid g\in \sym\rangle$, the classification of the $\sym$-pairs implies that  \emph{almost} everything is known 
about the semigroups 
$\langle g^{-1}ag\mid g\in G\rangle$, for each $\sym$-pair $(a,G)$ (see Section \ref{illustration}).

%

Regarding semigroups, our main result is the following.

\begin{thm}\label{sympairs}
Consider the pair  $(a,G)$, where $a$ is rank $r$ transformation with kernel partition of
type $\lambda$, and $G\leq \sym$. Then $(a,G)$ is an $\sym$-pair if and
only if one of the following holds:
\begin{enumerate}
\item $G$ is $\sym$ or $\alt$;
\item $r=1$ and $G$ is transitive;
\item $r>n/2$, $\lambda=(n-r+1,1,\ldots,1)$, and $G$ is $(n-r)$-homogeneous;
\item $r=n-2$, $\lambda=(2,2,1,\ldots,1)$, and $G$ is $4$-transitive;
\item $r=n-3$, $\lambda=(3,2,1,\ldots,1)$, and $G$ is $5$-transitive;
\item\label{(d)} $\lambda=(n-t,k_1,\ldots,k_{r-1})$ with $t<n/2$, $G$ is
$t$-homogeneous, and $(G,\lambda)$ is standard (as defined before Theorem
\ref{lambdatrans});
\item \label{(e)} finally we have a number of exceptional groups of low degree:\begin{enumerate}
\item $n=5$, $G$ is $\agl(1,5)$ and $\lambda\neq (2,2,1)$;
\item $n=6$, $G$ is $\psl(2,5)$ and $\lambda \not\in\{ (3,2,1),(3,1,1,1),(3,2,2),(2,2,1,1)\}$;
\item $n=6$, $G$ is $\pgl(2,5)$ and $\lambda \not\in\{ (2,2,1,1),(2,2,2)\}$;
\item $n=9$, $G$ is $\pgaml(2,8)$ and $\lambda$ {\bf is not one} of the following types:
\[\begin{array}{lllll}
( 2, 2, 2, 1, 1, 1 )& ( 2, 2, 2, 2, 1 )& ( 3, 2, 1, 1, 1, 1 )& ( 3, 2, 2, 1, 1 )&\\ 
 ( 3, 2, 2, 2 )& ( 3, 3, 1, 1, 1 )& ( 3, 3, 2, 1 )& ( 3, 3, 3 )&\\
   ( 4, 2, 1, 1, 1 )& ( 4, 2, 2, 1 )& ( 4, 3, 1, 1 )& ( 4, 3, 2 )& ( 4, 4, 1 ).
\end{array}
\]
\item $n=9$, $G$ is $\pgl(2,8)$ and $\lambda$ {\bf is one} of the following types:
\[\begin{array}{lllll}
( 2, 1, 1, 1, 1, 1, 1, 1 )& ( 3, 1, 1, 1, 1, 1, 1 )& ( 4, 1, 1, 1, 1, 1 )& ( 5, 1, 1, 1, 1 )\\
 ( 5, 3, 1 )& ( 5, 4 )& ( 6, 1, 1, 1 )& ( 6, 2, 1 )\\
  ( 6, 3 )& ( 7, 1, 1 )& ( 7, 2 )& ( 8, 1 ).
\end{array}
\]
\end{enumerate}
\end{enumerate}
\end{thm}

\medskip

We now outline the content of the paper.
 
In Section \ref{handt} we recall some well known results on homogeneity and transitivity. In Section \ref{prelim} we generalize  the well known correspondence between orbitals and orbits of a point stabilizer. We also give some character theory that is necessary later. Section \ref{hom} is devoted to the proof of Theorem \ref{lambdatrans}. A characterization of the $\sym$-pairs is proved in Section \ref{char}. Finally, Theorem \ref{sympairs} is proved in Section \ref{large}. The paper concludes with some remarks on pairs $(a,G)$ with $a$ a permutation
and $\langle a,G\rangle=\sym$ (Section \ref{perm}), with an illustration of the results on semigroups that immediately follow from the classification of the $\sym$-pairs (Section \ref{illustration}), and with several open problems in Section \ref{prob}.

\section{Homogeneity and transitivity}\label{handt}

We recall here the ``classical'' results on homogeneity and transitivity, since
we will make frequent use of them. Let $G$ be a subgroup of $\sym$.
Note that, if $\lambda=(n-t,1,\ldots,1)$ and $\mu=(n-t,t)$ then
\begin{itemize}
\item $G$ is $t$-homogeneous if and only if it is $\lambda$-homogeneous, and
$t$-transitive if and only if it is $\lambda$-transitive;
\item $G$ is $t$-homogeneous if and only if it is $\mu$-transitive (and this
is equivalent to $\mu$-homogeneity if $n\ne2t$).
\end{itemize}
In particular, $t$-homogeneity and $(n-t)$-homogeneity are equivalent; so, for
the purposes of classification, we may assume that $t\le n/2$.

A permutation group $G$ of degree $n$ is \emph{set-transitive} if it is
$t$-homogeneous for all $t$ with $1\le t\le n-1$. The problem of classification
of set-transitive groups was first posed by von Neumann and Morgenstern
\cite{vNM} in the context of fair $n$-player games, and was solved by
Beaumont and Petersen \cite{bp}: apart from the symmetric and alternating
groups there are just four set-transitive groups, namely
\begin{itemize}
\item $n=5$, $G=\agl(1,5)$;
\item $n=6$, $G=\pgl(2,5)$;
\item $n=9$, $G=\pgl(2,8)$ or $\pgaml(2,8)$.
\end{itemize}
Note that, if $G$ is set-transitive of odd degree $2k-1$, and we adjoin a new
point fixed by $G$, we obtain a $(k,k)$-transitive group. Three of the examples
in Theorem \ref{lambdahom} arise in this way.

Now $t$-homogeneity is equivalent to $(n-t)$-homogeneity, so we may assume
that $t\le n/2$. Theorem 1 of Livingstone and Wagner \cite{lw} shows
that $G$ is then $(t-1)$-homogeneous. Hence, if $G$ is not set-transitive,
then there is a unique integer $t<n/2$ (with possibly $t=0$) such that $G$ is
$s$-homogeneous if and only if either $s\le t$ or $s\ge n-t$. We say that
such a group is \emph{exactly $t$-homogeneous}.

Theorem 2(a) of Livingstone and Wagner \cite{lw} shows that a $t$-homogeneous group for
$t\le n/2$ is $(t-1)$-transitive. Therefore, it is immediate to see  that a
permutation group $G\leq \sym$ falls into one of the following classes:
\begin{enumerate}
\item $G$ is set-transitive, or 
\item $G$ is exactly $t$-homogeneous, for some $t\leq n/2$; in such  case:
\begin{enumerate}
\item either $G$ is $t$-transitive, but not $(t+1)$-homogeneous; or
\item $G$ is $t$-homogeneous, but not $t$-transitive.
\end{enumerate}
\end{enumerate}

Now the classification of $t$-homogeneous groups can be summarized as follows:
\begin{itemize}
\item if $5\le t\le n/2$ and $G$ is $t$-homogeneous, then $G$ is $t$-transitive
(Theorem 2(b) of Livingstone and Wagner~\cite{lw});
\item if $G$ is $2$-homogeneous but not $2$-transitive, then $G\le\agaml(1,q)$,
where $q$ is a prime power congruent to $3\pmod{4}$ (Kantor~\cite{kantor2});
\item if $G$ is $3$-homogeneous but not $3$-transitive, then either
$G\le\pgaml(2,q)$, where $q$ is a prime power congruent to $3\pmod{4}$, or
$G=\agl(1,8)$, $\agaml(1,8)$ or $\agaml(1,32)$, with $n=8,8,32$ respectively
(Kantor~\cite{kantor2});
\item if $G$ is $4$-homogeneous but not $4$-transitive, then $G=\pgl(2,8)$,
$\pgaml(2,8)$ or $\pgaml(2,32)$, with $n=9,9,33$ respectively
(Kantor~\cite{kantor});
\item using the Classification of Finite Simple Groups, all the $t$-transitive
groups for $t\ge2$ have been determined; lists are available in
\cite{cam,dixon}. In particular, the only $6$-transitive groups are the
symmetric and alternating groups, and the only $4$- or $5$-transitive groups
are these and the Mathieu groups $M_{11}$, $M_{12}$, $M_{23}$ and $M_{24}$.
\end{itemize}

\section{Preliminaries}\label{prelim}

The aim of this section is to prove  a result that generalizes the well known correspondence between orbitals and orbits of a point stabilizer. We also give a
character-theoretic interpretation. This material is well known to some group theorists but possibly not to semigroup theorists.

\begin{lem}\label{hig}
Let $A$ and $B$ be two sets and let $G$ be a group acting transitively on both
$A$ and $B$. Given $a\in A$, denote by $G(a)$ the stabilizer of $a$ in $G$. 

Then there is a natural bijection between the orbits of $G$ on $A\times B$ and the orbits in the action of $G(a)$ on $B$. In particular, the numbers of orbits
are equal. 
\end{lem}

\begin{pf}
Let $O$ be an orbit of $G$ on $A\times B$. For each $a\in A$, let $O(a)$
denote the set $\{b\in B:(a,b)\in O\}$.

We claim that $O(a)$ is an orbit of $G(a)$. For, if $b_1,b_2\in O(a)$, then
there exists $g\in G$ with $(a,b_1)g=(a,b_2)$; then $ag=a$ (so $g\in G(a)$)
and $b_1g=b_2$. Conversely, if $g\in G(a)$ and $b_1g=b_2$,
then $(a,b_1)g=(a,b_2)$.

Every $G(a)$-orbit in $B$ arises in this way: for, if $\Delta$ is such an
orbit, $b\in\Delta$, and $O$ is the $G$-orbit containing $(a,b)$, then
$\Delta=O(a)$. This completes the proof.
\qed
\end{pf}

\begin{figure}[h]
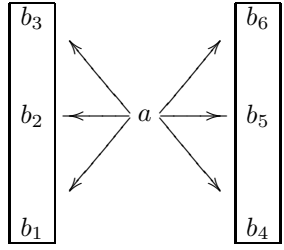

\[
\xy
(10,10)*{a}="a";
(12,9.6)*{}="a1";
(13,10)*{}="a2";
(12,10.4)*{}="a3";
(8,9.6)*{}="a4";
(7,10)*{}="a5";
(8,10.4)*{}="a6";
(20,0)*{}="b1";
(20,10)*{}="b2";
(20,20)*{}="b3";
(0,0)*{}="b4";
(0,10)*{}="b5";
(0,20)*{}="b6";
"a1";"b1" **\crv{} ?>* \dir{>};
"a2";"b2" **\crv{} ?>* \dir{>};
"a3";"b3" **\crv{} ?>* \dir{>};
"a4";"b4" **\crv{} ?>* \dir{>};
"a5";"b5" **\crv{} ?>* \dir{>};
"a6";"b6" **\crv{} ?>* \dir{>};
(-5,-5)*{b_{1}};
(-5,10)*{b_{2}};
(-5,23)*{b_{3}};
(25,-5)*{b_{4}};
(25,10)*{b_{5}};
(25,23)*{b_{6}};
(-8,-7)*{}="A"; (-8,25)*{}="B"; (-2,25)*{}="C";(-2,-7)*{}="D";
"A";"B" **\dir{-};
"B";"C" **\dir{-};
"C";"D" **\dir{-};
"A";"D" **\dir{-};
(22,-7)*{}="A1"; (22,25)*{}="B1"; (28,25)*{}="C1";(28,-7)*{}="D1";
"A1";"B1" **\dir{-};
"B1";"C1" **\dir{-};
"C1";"D1" **\dir{-};
"A1";"D1" **\dir{-};
\endxy
\]
\caption{A local graph of the action of $G$ on $A\times B$}\label{f2}
\end{figure}

One way of looking into this result is shown in the figure below. Represent $(a,b)\in A\times B$ by $a\rightarrow b$, and put a square around $b$ and $b'$ in case $(a,b)$ and $(a,b')$ are in the same orbit under $G$; then, locally, that is representing only the arrows rooted on $a$, we get figure \ref{f2}.

The fact that the local graph rooted on $a$ has  two orbits, means that the whole action of $G$ on $A\times B$ has two orbits. (Observe that the whole graph is obtained by translations of this picture under $G$ and hence the number of orbits is not going to change).  And looking at the picture is self-evident that $(a,b_{2})$ is in the orbit of $(a,b_{1})$ (under $G$) if and only if  $b_{2}$ is in the orbit of $b_{1}$ under $G(a)$. Thus the rectangles identify both the orbits under $G$ and under $G(a)$.

\bigskip

We can also use character theory to look at this result. The orbit-counting lemma (OCL for short), which states that the number of
orbits of a finite group $G$ on a finite set $A$ is equal to the average
number of fixed points of elements of the group, has a character-theoretic
interpretation. The \emph{permutation character} $\pi_A$ is the function on $G$
giving the numbers of fixed points of elements. OCL asserts that the number
of orbits of $G$ on $A$ is equal to the inner product
$\langle\pi_A,1_G\rangle_G$ of $\pi_A$ with the principal character $1_G$ of
$G$; in other words, the multiplicity of $1_G$ as a constituent of $\pi_A$.

We remark (this will be needed later) that the number of orbits of $G$ on
$A\times B$ is equal to the inner product of the permutation characters
of $G$ on the sets $A$ and $B$. For the permutation character on $A\times B$
is $\pi_A\pi_B$, and so
\[\begin{array}{rcll}
\hbox{\# $G$-orbits in $A\times B$} &=& \langle\pi_A\pi_B,1_G\rangle_G & \hbox{ (OCL)}\\
&=& \langle\pi_A,\pi_B\rangle_G & \hbox{ ($\pi_B$ is real).}
\end{array}\]
Note also that the proposition can be interpreted in terms of Frobenius
reciprocity. For the transitivity of $G$ on $A$ means that $\pi_A=1_{G(a)}|^G$
(the character of $G$ induced from the principal character of $G(a)$); so we 
have 
\[\begin{array}{rcll}
\hbox{\# $G(a)$-orbits in $B$} &=& \langle\pi_B|_{G(a)},1_{G(a)}\rangle_{G(a)} &
\hbox{ (OCL)}\\
&=& \langle\pi_B,1_{G(a)}|^G\rangle_G & \hbox{ (Frobenius reciprocity)}\\
&=& \langle\pi_B,\pi_A\rangle_G.
\end{array}\]

\begin{cor}
Suppose that the group $G$ acts transitively on two sets $A$ and $B$. Then
the stabilizer of a point in $A$ is transitive on $B$ if and only if the
stabilizer of a point in $B$ is transitive on $A$.
\label{thegreat}
\end{cor}

\begin{pf}
Both conditions are equivalent to transitivity of $G$ on $A\times B$, on applying Lemma \ref{hig} twice.
\qed
\end{pf}

This has a number of consequences which are relevant to our investigation, for
example the following.

\begin{cor}
Let $G\leq \sym$ be a $2$-homogeneous group such that its set-stabilizer of a $2$-set is transitive on the remaining points. Then the stabilizer of a point $p$ is transitive on the $2$-sets not containing $p$.  
\end{cor}

\section{Partition-homogeneous groups}\label{hom}

In this section we prove Theorem~\ref{lambdahom}, the classification of
$\lambda$-homogeneous groups which are not $\lambda$-transitive, and conclude
with some remarks about the classification of $\lambda$-transitive groups.

Until the end of the section, $\lambda$ denotes a partition of $n$ which is
not $(n)$ or $(1,1,\ldots,1)$, and $G$ a subgroup of $\sym$ which is
$\lambda$-homogeneous.

\subsection{Intransitive groups}\label{intransitive}

Suppose that $G$ is $\lambda$-homogeneous but intransitive. Then $G$ has at least two different orbits, and hence one of them must have at most $n/2$ elements. Therefore $G$ fixes
a set of size $t$, where $t\le n/2$. Let $T$ be the largest set fixed by $G$,
and such that $|T|=t\leq n/2$. So the stabiliser of $T$ in
$\sym$ is $\lambda$-homogeneous (because the $\lambda$-homogeneous group $G$ is contained in that stabilizer, and any over-group of a $\lambda$-homogeneous group is $\lambda$-homogeneous). Therefore  the stabiliser $S$ of a
partition of shape $\lambda$ is $t$-homogeneous (by Corollary \ref{thegreat}). The stabiliser of a
partition is (by definition) imprimitive. Now suppose, by contradiction, that $t\geq 2$; then, as $S$ is $t$-homogeneous, by Theorem 1 
of Livingstone and Wagner \cite{lw}, it follows that $S$ is $2$-homogeneous, and hence primitive. 
So necessarily $t=1$. By choice of $t$,  we see that $G$ is transitive on the
remaining points.

In this case, the partition $\lambda$ must be uniform (all parts have the same
size); for, if not, then two set partitions in which the fixed point lies in
parts of different sizes could not be equivalent under $G$. 
We will return to this case later.

\subsection{Uniform partitions}\label{uni}

The aim of this subsection is to prove the following theorem.
\begin{thm}\label{lem1}
Let $\lambda$ be a uniform partition of $n$ other than $(1,1,\ldots,1)$, 
and let $G$ be a $\lambda$-homogeneous group.
Then one of the following holds:
\begin{enumerate}
\item $G$ is $\lambda$-transitive;
\item $n=6$, $\lambda =(3,3)$ and either $G$ is $\psl(2,5)$, or $G=\agl(1,5)$ (fixing a point);
\item $n=10$, $\lambda=(5,5)$ and $G$ is $\psl(2,8)$ or $\pgaml(2,8)$ (in both cases fixing a point);
\item  $G$ fixes a point and acts as the alternating or the symmetric group on the remaining points.
\end{enumerate}

\end{thm}

The proof of this theorem will be carried out in a number of paragraphs.  
As said above,  $\lambda$ denotes a uniform partition, say
$\lambda=(k,k,\ldots,k)$ ($t$ parts), with $n=tk$, and $t,k>1$.

\textbf{The groups of degree at most 4}

By direct inspection of the groups of degree up to $4$ we see that the only $\lambda$-homogeneous groups (when $\lambda$ is uniform) are $\sym$ and $\alt$.  

\textbf{Transitivity}

The aim of this paragraph is to prove that either $G$ is transitive or it fixes one point $z$ and is $2$-homogeneous on $\Omega\setminus\{z\}$.

We saw above that, if $G$ is intransitive, then it fixes a point $z$ and is
transitive on the remaining points. Each partition has a unique part
containing $z$, of the form $\{z\}\cup L$, where $|L|=k-1$. So $G$ must be
$(k-1)$-homogeneous on $\Omega\setminus\{z\}$. If $k>2$, then it follows that
$G$ is $2$-homogeneous on $\Omega \setminus \{z\}$.

In the case $k=2$, we need a little more. We can assume that $t>2$, since the
 groups of degree~$4$ are treated above. The stabiliser of a
$2$-set $\{z,x\}$ is $\mu$-homogeneous on the remaining points, where $\mu$ has 
$t-1$ parts of size~$2$. If this stabiliser $H$ is not transitive, then (by
the preceding argument) it fixes a point $x'$; so $\{x,x'\}$ is a block of
imprimitivity for $G$ acting on $\Omega\setminus\{z\}$. But this is impossible,
since $2$ does not divide $n-1=2t-1$. So the stabiliser of $x$ is transitive
on the remaining points, and $G$ is $2$-transitive on $\Omega\setminus\{z\}$.

Hence, under the conditions of Theorem \ref{lem1}, either $G$ is transitive,
or it fixes one point $z$ and is $2$-homogeneous on $\Omega\setminus\{z\}$.

\textbf{$2$-homogeneity} Next we show that, if $G$ is transitive, then it
is $2$-homogeneous. This involves a little character theory, similar to (but
much simpler than) that used by Livingstone and Wagner~\cite{lw} and Martin
and Sagan~\cite{ms}; see the comments before Corollary~\ref{thegreat}.

Consider the action of the symmetric group on $\Omega$, the set of $2$-subsets
of $\Omega$, and the set of partitions of $\Omega$ of type $(k,k,\ldots,k)$.
It is well known that the permutation characters in the first two of these
actions have the form $\chi_0+\chi_1$ and $\chi_0+\chi_1+\chi_2$, where
$\chi_0$ is the principal character, and $\chi_1$ and $\chi_2$ are the
irreducibles corresponding to the partitions $(n-1,1)$ and $(n-2,2)$
respectively. The stabiliser of a partition is transitive on points and has
two orbits on $2$-sets, so the permutation character $\pi$ of the action on
partitions contains $\chi_0$ and $\chi_2$ (but not $\chi_1$).

On restriction to $G$, the character $\pi$ contains the principal character
with multiplicity~$1$ (since $G$ is transitive on the set of partitions). So
the principal character is not contained in $\chi_2|_G$. Also it is not
contained in $\chi_1|_G$, since $G$ is assumed transitive. So $G$ is
transitive in its action on $2$-sets, that is, $2$-homogeneous.

\textbf{Completion of the proof of Theorem \ref{lem1}}
To complete the proof we make use of three facts. The first two are
versions of Bertrand's postulate, which are far from the strongest form known;
an accessible proof of what we need can be found on Robin Chapman's homepage
\cite{rjc}. The third is a theorem of Jordan (see Wielandt~\cite{wie}, Theorem 13.9).
\begin{enumerate}
\item If $n\ge8$, then there is a prime $p$ satisfying $n/2<p\le n-3$.
\item If $n\ge16$, then there is a prime $p$ satisfying $n/2<p\le n-4$.
\item If $G$ is a primitive permutation group of degree $n$ containing a
$p$-cycle, where $p$ is a prime satisfying $p\le n-3$, then $G$ is the
symmetric or alternating group.
\end{enumerate}

The number of partitions of $\{1,2,\ldots,n\}$ of type $(k,k,\ldots,k)$, with
$n=tk$,  is equal to $\frac{n!}{(k!)^t\ t!}$,  which is divisible by a prime
$p$ satisfying the conditions of
(1) above, provided that $n\ge8$. So $|G|$ is divisible by $p$, and $G$
contains an element of order $p$, necessarily a $p$-cycle (as $p>\frac{n}{2}$). 

Suppose first that
$G$ is transitive and $n\geq 8$. Then it is $2$-homogeneous, and hence primitive, and so by
Jordan's theorem it is symmetric or alternating. 

For $4<n<8$, $G$ transitive,
and $\lambda$ uniform, the only possible degree is $n=6$ (as the other are primes). Therefore the only possible groups are $\psl(2,5)$ and $\pgl(2,5)$. The latter is $3$-homogeneous and hence $\lambda$-transitive for $\lambda=(3,3)$. The former is $2$-homogeneous   
and this yields the exception $G=\psl(2,5)$, for $\lambda =(3,3)$. 

If $G$ is intransitive, then it fixes
a point and is $2$-homogeneous (and so primitive) on the remaining $n-1$ points;
so, if $n\ge16$, we conclude that $G$ is the symmetric or alternating group on
the moved  points. 

To finish the proof of the theorem we only need to inspect the intransitive
groups of degree $n$ for $n\le16$. If such a group $G$ is $\lambda$-homogeneous
for $\lambda=(k,k,\ldots,k)$ ($t$ parts), then the number $\frac{n!}{(k!)^t\ t!}$
of set partitions of shape $\lambda$ must divide $|G|$. Moreover, by the results in Section \ref{handt}, we know that a $2$-homogeneous group 
 is either $2$-transitive or primitive of odd order. A quick check of
the primitive groups using 
GAP \cite{GAP} shows that the three groups of the theorem are the only possibilities. The theorem is proved.

Note that, in cases (2)--(4) of the theorem, the partition is $\lambda=(k,k)$
with $n=2k$, and so $\lambda$-homogeneity is equivalent to $(k-1)$-homogeneity,
that is, to set-transitivity, on the non-fixed points.

\subsection{The non-uniform case}

In this subsection we will prove Theorem \ref{lambdahom} for non-uniform
partitions. From Subsection \ref{intransitive} we know that if $G$ is $\lambda$-homogeneous and intransitive, then $\lambda$ is uniform; therefore, a $\lambda$-homogeneous group, for a non-uniform $\lambda$, must be transitive. In fact more can be said.

\begin{lem}\label{imprimiti}
Let $\lambda$ be a non-uniform  partition of $n$, and let $G$ be a
$\lambda$-homogeneous but not $\lambda$-transitive group. Then $G$ is
$2$-homogeneous, and hence primitive.
\end{lem}

\begin{pf}
Let $\lambda'$ be the partition obtained from $\lambda$ by summing all
the occurrences of each number. Then $G$ is $\lambda'$-transitive, from
which it follows that, if $r$ is the sum of any subset of the parts of
$\lambda$, then $G$ is $r$-homogeneous. If there is such an $r$ with
$1<r<n-1$, then $G$ is $2$-homogeneous, as required; so suppose not. Then
$\lambda'=(n-1,1)$, so that $\lambda=(k,\ldots,k,1)$ (with $t$ parts
equal to $k$, with $tk=n-1$). As $G$ is not $\lambda$-transitive, it follows that $\lambda\ne(n-1,1)$, so $1<k,t<n-1$.

Now $G(x)$, the stabiliser in $G$ of a point $x$,  is $(k,\ldots,k)$-homogeneous on the
remaining points. If it is transitive, then $G$ is $2$-transitive. Otherwise
it fixes a point $z$ and is $2$-homogeneous (and hence primitive) on the
points different from $x$ and $z$. Now $\{x,z\}$ is a block of imprimitivity
for $G$; the primitivity of $G(x)$ implies that there are only two blocks,
so $n=4$, which is impossible. The lemma follows.
\qed
\end{pf}

\vspace{0.2cm}

Now we prove Theorem \ref{lambdahom} (for non-uniform partitions).  Let $G$ be $\lambda$-homogeneous, where $\lambda$ has $r_k$ parts equal
to $k$, for $k=1,2,\ldots$. Then the order of $G$ is divisible by the number
$\frac{n!}{\prod_k(k!)^{r_k}\,r_k!}$ of unordered set partitions of
shape $\lambda$.

Suppose that $n\ge8$. By Bertrand's Postulate, there is a prime $p$
satisfying $n/2<p\le n-3$. If $p\mid|G|$, then by Jordan's Theorem, $G$ is
symmetric or alternating, and so is $\lambda$-transitive. So we may assume 
that $p\nmid|G|$, whence also $p$ does not divide the number of set partitions.
This implies that either $k\ge p$ or $r_k\ge p$ for some $k$. Since
$kr_k\le n$ and $p>n/2$, we have either $k\ge p$, $r_k=1$, or $k=1$,
$r_k\ge p$; in other words, either $\lambda$ has a unique part of size at
least $p$, or it has at least $p$ parts equal to $1$.

Denote by $t$ the number such that, in the above, either $k=n-t$ or
$r_k=n-t$. Then $G$ is $t$-homogeneous, so $t\le 5$. The group $H$ induced
on a $t$-set by its setwise stabiliser is $\lambda'$-transitive, where
$\lambda'$ is obtained from $\lambda$ by removing either the part $n-t$ or
the $n-t$ parts $1$.

Suppose first that $\lambda$ has $n-t$ parts equal to $1$. Then $\lambda'$
has no part equal to $1$; so $\lambda'=(3,2)$ or $(2,2)$. In the first case,\
$G$ is $5$-homogeneous, so is $M_{12}$ or $M_{24}$; in the second case, $G$
is $4$-homogeneous, and $H$ is $(2,2)$-homogeneous or $(2,1,1)$-homogeneous,
so has order divisible by $3$; thus $G$ is $M_n$ (for $n=11, 12, 23, 24$)
or $\pgaml(2,8)$. All these groups occur in the statement of the theorem.

Now suppose that $\lambda$ has a part equal to $n-t$. Then $G$ is
$t$-homogeneous but not $t$-transitive, so $t\le 4$; and $H$ is
$\lambda'$-homogeneous but not $\lambda'$-transitive. In the case
$\lambda'=(1,\ldots,1)$, any $t$-homogeneous but not $t$-transitive group
can occur. Otherwise, since the parts of $\lambda'$ cannot all be distinct,
we have $\lambda'=(2,2)$ or $(2,1,1)$. These cases occur only for the group
$G=\pgaml(2,8)$, which is $\lambda$-transitive.

Finally, the cases where $n\le7$ are easily checked by hand or using GAP.

The proof of Theorem \ref{lambdahom} is finished.

\subsection{$\lambda$-transitive groups}\label{eti}

In this section we prove Theorem~\ref{lambdatrans}. So $G$ is a subgroup
of $\sym$, not $\sym$ or $\alt$, which is $\lambda$-transitive, where
$\lambda\ne(n)$.

As observed before the theorem, it is clear that, if $(G,\lambda)$ is standard, then $G$ is
$\lambda$-transitive; and conversely, if the largest part of $\lambda$ is
$n-t$, with $t\le n/2$, and $G$ is $\lambda$-transitive, then $G$ is standard. Therefore the only thing still required to prove Theorem~\ref{lambdatrans} is that, if $G$ is $\lambda$-transitive, then the largest part of $\lambda$ must be greater than $n/2$. 

Again, if $n\ge8$, there is a prime $p$ satisfying $n/2<p\le n-3$. If
$p\mid|G|$, then $G$ is symmetric or alternating, contradicting our assumptions; we may assume that this
is not the case. So $p$ does not divide the number of \emph{ordered}
set partitions of type $\lambda$, which is $n!/\left(\prod_ik_i!\right)$,
where $\lambda=(k_1,k_2,\ldots)$. So $k_1\ge p>n/2$. Our comments
above show that the pair $(G,\lambda)$ is standard in this case.

Finally, the groups of degree at most $7$ are easily handled directly.

\section{A characterization of the $\sym$-pairs}\label{char}

 The goal of this section is to provide a characterization of the $\sym$-pairs, that is, 
 the pairs $(a,G)$, where $a\in \trans\setminus G$ and $G\leq \sym$, such that 
\[
\langle a,G\rangle \setminus = \langle a,\sym\rangle \setminus \sym .
\]

We start by introducing some notation and terminology. Recall that $\Omega=\{1,\ldots ,n\}$. The rank of a map $a\in \trans$ is the size of its image: $\rank(a)=|\Omega a|$; the  kernel of $a\in \trans$ is the partition of $\Omega$ induced by the equivalence relation $\{(x,y)\in \Omega\times \Omega\mid xa=ya\}$. The kernel type of $a$ is the tuple $(|A_{1}|,\ldots,|A_{k}|)$, where the sets $A_{i}$ belong to the kernel of $a$ and for $i\leq j$ we have $|A_{i}|\geq |A_{j}|$. 
\begin{lem}\label{lemma1}
Consider an $\sym$-pair  $(a,G)$, with $G\leq \sym$ and $a\in\trans$ such that $\ker(a)$ has type $\lambda$.  Then
\begin{enumerate}
\item $G$ is $\rank(a)$-homogeneous;
\item $G$ is $\lambda$-homogenous.
\end{enumerate}
\end{lem}
\begin{pf}
Let $a\in \trans\setminus\sym$. For every $A\subseteq \Omega$ of size $|\Omega a|$, there exists $g\in \sym$ such that $\Omega ag=A$ and hence there exists  $b\in\langle a,\sym\rangle\setminus \sym$ (namely $b=ag$) such that $\Omega b=A$. Therefore, as $(a,G)$ is an $\sym$-pair, it follows that for every $A\subseteq \Omega $ of size $|\Omega a|$, there exist $g_{1},\ldots, g_{k}\in G$ such that $\Omega g_{1}ag_{2}a\ldots ag_{k}=A$. As $\Omega g_{1}ag_{2}a\ldots a\subseteq \Omega a$ (and $|\Omega g_{1}ag_{2}a\ldots a|=|A|=|\Omega a|$ and we are working with finite sets) it follows that $\Omega g_{1}ag_{2}a\ldots a=\Omega a$.  Thus $g_{k}\in G$ maps $\Omega a$ onto $A$. It is proved that $G$ must be $|\Omega a|$-homogeneous.  

Similarly, let $P=(A_{1},\ldots,A_{m})$ be the partition of $\Omega$ induced by $\ker (a)$. Then given any partition $\pi$ of $\Omega$ with the same kernel type of $a$, there exists $g\in \sym$ such that  $Pg=\pi$. Therefore,  as $(a,G)$ is an $\sym$-pair, for every partition $\pi$ of $\Omega $, of the same type as $P$,  there exist $g_{1},\ldots, g_{k}\in G$ such that $\ker( g_{1}ag_{2}a\ldots ag_{k})=\pi$. Using an argument similar to the one used above one can see that $Pg^{-1}_{1}=\pi$.\qed   
\end{pf}
	
\smallskip
These two necessary conditions turn out to be sufficient as well.

\begin{lem} Let $a\in \trans$ such that $\rank(a)=k$ (with $k<n$) and kernel type $\lambda$.
Let $G\leq \sym$ be a  $k$-homogeneous and $\lambda$-homogenous group. Then $(a,G)$ is an $\sym$-pair.  
\end{lem}\label{lemma2}
\begin{pf}
Let $\Sigma_{a}=\langle \sym ,a\rangle \setminus \sym $ and  $\Gamma_{a}=\langle G,a \rangle\setminus G$.
Since $\Gamma_{a} \subseteq  \Sigma_{a}$, we only need to prove the converse. Observe that by theorems \ref{t3} and \ref{mainold} the semigroup $\Sigma_{a}$ is idempotent generated. In addition, it is obvious that  $\Sigma_{a}$ is generated by the elements $\{gah\mid g,h\in \sym\}$. As every element in this set has the same rank as $a$, it follows that the semigroup $\Sigma_{a}$ is generated by its idempotents  of the largest rank (those idempotents with the same rank as $a$). Hence, to prove the desired inclusion we only need to prove that all the idempotents of $\Sigma_{a}$ of the same rank as $a$ are contained in $\Gamma_{a}$.

Observe that given any pair $(\rho,S)$, where $\rho$ is a partition of $\Omega$ and $S$ is a section for $\rho$, there exists one and only one idempotent $e$ such that $\ker(e)=\rho$ and $\Omega e=S$. Therefore, if $e^{2}=e\in \Sigma_{a}$, with $\rank(e)=\rank(a)$, this means that  $e=g_{1}ag_{2}\ldots g_{k-1}ag_{k}$ ($g_{i}\in \sym)$; this implies that $\ker(e) \supseteq \ker (g_{1}a)$ and $\Omega e \subseteq \Omega ag_{k}$. As $\rank(e)=\rank(a)$ it follows that equalities hold so that $\ker(e)=\ker(g_{1}a)$ and $\Omega e= \Omega ag_{k}$. Using the well known fact that $\ker(ga)=\ker(a) g^{-1}$, we derive the following characterization of the idempotents possessing the same rank as $a$:   $e^{2}=e \in \Sigma_{a}$, with $\rank(a)=\rank(e)$, if and only if there exist $g,h\in\sym$ such that $\ker(e)=\ker(a)g$, $\Omega e = \Omega a h$, and $\Omega a h$ is a section for  $\ker(a)g$. Fix $g,h\in \sym$ such that $\Omega a h$ is a section for  $\ker(a)g$. To prove the theorem we only need to prove that  if $e$ is the unique idempotent associated to the pair $(\ker(a)g,\Omega a h)$, then $e\in \Gamma_{a}$. 

At this point a further simplification is possible. In fact, if it can be proved that there exists $b\in \Gamma_{a}$ such that $\ker(b)=\ker(a)g$ and $\Omega b= \Omega a h$, then any power of $b$ also belongs to $\Gamma_{a}$ (since $b\in \Gamma_{a }$ and $\Gamma_{a}$ is a semigroup). As the image of $b$ is a section for its kernel (as we are assuming that $\Omega b=\Omega a h$ is a section for  $\ker(b)=\ker(a)g$), it follows that $\rank(b^{k})=\rank(b)$, for all natural $k$, and hence ---by the result in \cite[pp. 9-11]{Ho95}---there exists a natural number (usually denoted by $\omega$) such that $b^{\omega}$ is an idempotent with kernel equal to $\ker (b)$ and with image equal to $\Omega b$. Therefore, to prove that the unique idempotent $e$ such that $\ker(e)=\ker(a)g$ and $\Omega e=\Omega a h$ belongs to $\Gamma_{a}$, it is sufficient to prove that there exists $b\in \Gamma_{a}$ such that $\ker (b)=\ker(a) g$ and $\Omega b=\Omega a h$. But now the conclusion follows immediately. Since $G$ is $|\Omega a|$-homogeneous, there exists $h'\in G$ such that $\Omega a h'=\Omega a h$; similarly, as $G$ is $\lambda$-homogenous, there exists $g'\in G$ such that $\ker(a)g'=\ker(a)g$. As $(g')^{-1}ah'\in \Gamma_{a}$ it follows that $b=(g')^{-1}ah'$ has the desired properties: $\ker(b)=\ker(a)g'=\ker(a)g$ and $\Omega b=\Omega a h' =\Omega a h$. 
\qed
\end{pf}

\smallskip

The two previous lemmas yield the following result which is the main theorem in this section.

\begin{thm}\label{main2}
Let $a\in\trans$ such that $\rank(a)=k$ and $\ker(a)$ has type $\lambda$. Then  $(a,G)$, with $G\leq \sym$, is an $\sym$-pair if and only if $G$ satisfies the following two properties: 
\begin{enumerate}
\item $G$ is $k$-homogeneous; and 
\item $G$ is $\lambda$-homogenous.
\end{enumerate}
\end{thm}

\section{The classification of the $S_{n}$-pairs}\label{large}

The aim of this section is to prove the classification theorem of the $\sym$-pairs (Theorem~\ref{sympairs}).

That $(a,\sym)$ is an $\sym$-pair for every non-invertible transformation is trivial. That  $(a,\alt)$ is also an $\sym$-pair for every non-invertible transformation was proved by Levi \cite{levi96}. It is very easy to see that,  for a constant map $a$,   $(a,G)$  is an $\sym$-pair if and only if $G$ is transitive. 

Suppose that $a\in\trans\setminus\sym$, with  $\rank (a)=r>1$ and kernel partition of
type $\lambda$, and that $(a,G)$ is an $\sym$-pair. Then $G$ is
$r$-homogeneous and $\lambda$-homogeneous by Theorem \ref{main2}. In particular, $G$ is transitive.
We have to look at the groups in the first two theorems.

Suppose $G$ is $r$-homogeneous and $\lambda$-homogeneous, but not $\lambda$-transitive; then $G$ satisfies one of the conditions of  Theorem \ref{lambdahom} (\ref{trans}). Suppose $\lambda=(n-t,1,\ldots ,1)$. As $\lambda$ is the kernel type of a $\rank r$ map, this means that $n-t=n-r+1$, that is, $t=r-1$. Since, by    Theorem \ref{lambdahom} (\ref{trans}), $G$ is $t$-homogeneous, but not $t$-transitive,  it follows that $G$ is $(r-1)$-homogeneous, but not $(r-1)$-transitive. However, $G$ is $r$-homogeneous; thus $r>n/2$.


The other
groups of the theorem are among the required conclusions.

Suppose now that $G$ is $\lambda$-transitive, so that $(G,\lambda)$ is standard.
This gives conclusion (\ref{(d)}) of the theorem.

The small groups (the set-transitive groups together with $\psl(2,5)$) occur
under conclusion (\ref{(e)}).

\section{Permutations}\label{perm}

In the non invertible case it does not make sense to ask for the transformations $a\in \trans$ such that $\langle a,G\rangle=\langle a,\sym\rangle$ since this implies $G=\sym$. Therefore we had to search for the pairs  $(a,G)$ such that $\langle a,G\rangle\setminus G=\langle a,\sym\rangle\setminus \sym$.

Regarding permutations $a\in \sym$, it does not make sense to ask for the groups $G$  such that $\langle a,G\rangle \setminus G = \langle a,\sym\rangle \setminus \sym$ since this happens if and only if $a\in G$. In the permutations setting the analogous of the non invertible problem is to ask for pairs $(a,G)$ with
$\langle a,G\rangle=\langle a,\sym\rangle=\sym$. A great deal is known about this problem.

If $G$ is the trivial group, then there is no such permutation $a$, provided
$n>2$, since the symmetric group requires at least two generators. However,
a theorem of {\L}uczak and Pyber~\cite{lp} shows that, if $G$ is transitive,
then almost all $a\in\sym$ satisfy $\langle a,G\rangle=\sym$ or $\alt$
(in the sense that the proportion of $a$ satisfying this tends to $1$ as
$n\to\infty$). Of course, if $\langle a,G\rangle=\alt$, then $a$ is an
even permutation and $G$ contains only even permutations.

It is interesting to note that this theorem, together with a
theorem of H\"aggkvist and Janssen~\cite{hj}, answers a question about
quasigroups: it shows that almost all quasigroups of order $n$ have the
property that the right multiplication group is $\sym$. 

There is an extensive literature about the problem in the case where $G$ is
a cyclic group, which we do not survey here.

\section{Consequences of the classification of the $\sym$-pairs}\label{illustration}

As observed above much is known  on the $\sym$-normal semigroups and hence the characterization of the $\sym$-pairs gives for free  a large number of results. The aim of this section is to have a glance on those consequences. We start with a number of easy observations.

We recall the following result \cite{symo}. 
\begin{thm}\label{symons}
Let $a\in\trans \setminus \sym$. Then 
\[
\langle a,\sym \rangle\setminus \sym = \{b\in \trans \mid (\exists g\in \sym)\ \ker(a)g\subseteq \ker(b)\}.
\]
\end{thm}

As we have  the classification of $\sym$-pairs, now we would like to investigate what happens with sets $A\subseteq \trans\setminus \sym$. It would be good to have a result of the form $\langle A,G\rangle\setminus G=\langle A,\sym\rangle\setminus \sym$  if and only if $(a,G)$ is an $\sym$-pair, for all $a\in A$. However this is not true. 
 The smallest example, by Theorem \ref{sympairs}, is provided by $G:=\agl(1,5)$ of degree $5$, and $A=\{a,b\}$, where $a$ and $b$ are the following transformations: 

\[\begin{array}{l}
a=\left(
\begin{array}{cccccccc}
\{1,2\}&\{3\}&\{4\}&\{5\}\\
1&3&4&5
\end{array}\right)\\ \\
b=\left(
\begin{array}{cccccccc}
\{1,2\}&\{3,4\}&\{5\}\\
1&3&5
\end{array}\right).
\end{array}
\]

Since $(a,G)$ is an $S_{n}$-pair, it follows from Theorem \ref{symons}  that $\langle a,G\rangle\setminus G=\trans\setminus \sym$ and hence contains $b$; thus $\langle A,G\rangle \setminus G= \langle a,G\rangle \setminus G=\langle a,\sym\rangle \setminus \sym$.  But $(b,G)$, by Theorem \ref{sympairs}, is not an $S_{n}$-pair. 

The following result shows that regarding $\sym$-pairness only the top elements matter. 
\begin{lem}
Let $a\in \trans\setminus \sym$ and $G\le \sym$. The following are equivalent:
\begin{enumerate}
\item $\langle a,G\rangle\setminus G =\langle a,\sym \rangle \setminus \sym$;
\item $\{b\in \langle a,G\rangle \mid \rank(b)={\rank(a)}\}=\{b\in \langle a,\sym\rangle \mid \rank(b)={\rank(a)}\}$.
\end{enumerate}
\end{lem}
\begin{pf}
The forward direction is obvious. Regarding the converse, observe that, for any group $H$, the semigroup $\langle a,H\rangle\setminus H$ is generated by the set $\{hag\mid g,h\in H\}$, and this set is contained in $\{b\in \langle a,H\rangle \mid \rank(b)={\rank(a)}\}$. Therefore, $$\langle a,H\rangle\setminus H=\langle \{b\in \langle a,H\rangle \mid \rank(b)={\rank(a)}\}\rangle.$$

 Thus
\[
\begin{array}{rcl}
\langle a,G\rangle\setminus G&=&\langle  \{b\in \langle a,G\rangle \mid \rank(b)={\rank(a)}\} \rangle \\
&=&  \langle \{b\in \langle a,\sym\rangle \mid \rank(b)={\rank(a)}\} \rangle\\
& =& \langle a,\sym \rangle \setminus \sym.
\end{array}\] The result follows.~\qed
\end{pf}

\vspace{0.2cm}

A set $A\subseteq \trans\setminus \sym$ is said to be independent if for every $a,b\in A$ ($a\neq b$) we have $a\not \in \langle b,\sym\rangle$. Since, given any $A\subseteq \trans\setminus \sym$,  we can pick the maximal elements of $\{\langle a,\sym\rangle\setminus \sym \mid a\in A\}$, then every set $A\subseteq \trans\setminus \sym$ contains (at least) one independent subset. This set  is not unique just because we can have $a,b\in A$ such that  $\langle a,\sym\rangle\setminus \sym = \langle b,\sym\rangle\setminus \sym$. Theorem \ref{symons} immediately provides a characterization of independent sets. 

\begin{lem}\label{indep}
Let $A\subseteq \trans\setminus \sym$ be a set. The following are equivalent:
\begin{enumerate}
\item $A$ is independent;
\item for all different $a,b\in A$ there is no $g\in \sym$ such that $\ker(a)g\subseteq \ker(b)$.
\end{enumerate}
\end{lem}

The following theorem shows the nice behavior  of the idea of   $\sym$-pairs when extended to pairs $(A,G)$, for an independent set $A$. 

\begin{thm}
Let $A\subseteq \trans\setminus \sym$ be an independent set and let $G\le \sym$. The following are equivalent:
\begin{enumerate}
\item $\langle A\cup G\rangle \setminus G =\langle A\cup \sym \rangle \setminus \sym$;
\item   $(a,G)$ is an $\sym$-pair, for all $a\in A$. 
\end{enumerate}
\end{thm}
\begin{pf}
We start by proving the theorem in the forward direction. Let $a\in A$. If $a\in \langle (A\setminus \{a\})\cup G\rangle$, then $a=g_{1}a_{1}\ldots a_{k}g_{k+1}$ (with $a_{i}\in A\setminus \{a\}$, $g_{i}\in G$). This implies that  $\ker(a)\supseteq \ker(a_{1})g^{-1}_{1}$. Thus, by Theorem \ref{symons}, $a\in \langle a_{1},\sym\rangle $, a contradiction since $A$ is independent. 

Conversely, it is clear that $\langle A\cup G\rangle \setminus G \subseteq \langle A\cup \sym \rangle \setminus \sym$. So let $x\in \langle A\cup \sym \rangle \setminus \sym$. Then $x=g_{1}a_{1}g_{2}\ldots $, with $a_{1}\in A, g_{i}\in \sym$. Thus $\ker(x)\supseteq \ker(a_{1})g^{-1}_{1}$ and hence, by Theorem \ref{symons}, we have 
\[
x\in \langle a_{1},\sym\rangle \setminus \sym = \langle a_{1},G\rangle \setminus G\subseteq  \langle A,G\rangle \setminus G.  
\]
\qed
\end{pf}

We now state the main theorem of this illustrative section.

\begin{thm}\label{illustr}
Let $(a,G)$ be an $\sym$-pair and let $S=\langle a,G\rangle \setminus G$. Then 
\begin{enumerate}
\item\label{i1} $S= \{b\in \trans \mid (\exists g\in \sym)\ \ker(a)g\subseteq \ker(b)\}$;
\item\label{i2} $S$ is regular, that is, for all $a\in S$ there exists $b\in S$ such that $a=aba$;
\item\label{i3} $S$ is generated by its idempotents;
\item\label{i4} $S$ and $\langle g^{-1}ag \mid g\in G\rangle$ have the same idempotents;
\item\label{i5} $S=\langle g^{-1}ag \mid g\in G\rangle$;
\item\label{i6} $G$ synchronizes $a$ (that is $G$ and $a$ generate a constant transformation) and hence (as $G$ is transitive) $S$ contains all the constants;
\item\label{i7} since $S$ contains all the constants, the automorphisms of $S$ coincide with the elements of the normalizer of $S$ in $\sym$, $$\mbox{Aut}(S)=\mbox{N}_{\sym}(S).$$ 
\item\label{i8}  in fact we also have $\mbox{Aut}(S)\cong \sym$;
\item\label{i9} all the congruences of $S$ are described;
\item\label{i10} if $e^{2}=e\in S$, $r:=\rank(e)$, then $$\{f\in S \mid \ker(f)=\ker (e) \mbox{ and }\Omega f=\Omega e\}\cong \mathcal{S}_{r}.$$
\item\label{i11} regarding principal ideals and Green's relations, for all $a,b\in S$, we have \[
\begin{array}{rcl}
aS=bS &\Leftrightarrow & \ker(a)=\ker(b)\\	
Sa=Sb &\Leftrightarrow & \Omega a=\Omega b\\	
SaS=SbS &\Leftrightarrow & \rank (a)=\rank( b )\\
\end{array}
\]
\end{enumerate}
\end{thm}
\begin{pf}
 Equality (\ref{i1}) was proved by Symons in \cite{symo}. Claims (\ref{i2}), (\ref{i3}) and (\ref{i5}) were proved by Levi and McFadden in \cite{lm}. Claim (\ref{i4}) was proved by McAlister in \cite{mcalister}, and (together with (\ref{i3})) it also implies (\ref{i5}). 
 
 It is well known that every $2$-homogeneous group is synchronizing; since $\sym$ is $2$-homogeneous claim (\ref{i6}) follows. 
 
 Claim (\ref{i7}) belongs to the folklore, probably was first observed by Schreier \cite{Sc36}; see also \cite{arko} and \cite{arko2}.
 Claim (\ref{i8}) was proved by Symons in \cite{symo}, but is also an easy consequence from (\ref{i7}).
In \cite{levi00} Levi  described all the congruences of an $\sym$-normal semigroup  and hence described the congruences in $S$. Thus (\ref{i9}).

Stetement  (\ref{i10}) belongs to the folklore (see Theorem 5.1.4 of \cite{mar}).  
The results about principal ideals  (\ref{i11}) were proved by Levi and McFadden in \cite{lm}.

\qed
\end{pf}

\section{Problems}\label{prob}

We start with a problem about permutations. 
\begin{problem}
Let $H\leq \sym$ be a $2$-transitive group. Classify the pairs $(a,G)$, where $a\in \sym $ and $G\leq \sym$, such that  
$\langle a,G\rangle  = H$. 
\end{problem}

The main questions answered in this paper admit linear analogous. 

\begin{problem}
Let $V$ be a finite dimension vector space. A pair $(a,G)$, where $a$ is a singular endomorphism of $V$ and  $G\leq\Aut(V)$, is said to be an $\Aut(V)$-pair if $$\langle a,G\rangle\setminus G=\langle a,\Aut(V)\rangle\setminus\Aut(V).$$

Classify  the $\mbox{Aut}(V)$-pairs.
\end{problem}	

To handle this problem it is useful to keep in mind the following results. Kantor~\cite{kantor:inc} proved that if a subgroup of $\pgaml(d,q)$ acts transitively on $k$-dimensional subspaces, then it acts transitively on $l$-dimensional subspaces for all $l\le k$ such that $k+l\le n$; in~\cite{kantor:line}, he showed that subgroups transitive on $2$-dimensional subspaces are $2$-transitive on the $1$-dimensional subspaces with the single exception of a subgroup of $\pgl(5,2)$ of order $31\cdot5$; and, with the third author~\cite{cameron-kantor}, he showed that such groups must contain $\psl(d,q)$ with the single exception of the alternating group $A_7$ inside $\pgl(4,2)\cong A_8$. Also Hering \cite{He74,He85} and Liebeck \cite{Li86} classified the subgroups of $\pgl(d,p)$ which are transitive on $1$-spaces.

\vspace{0.2cm}

For definitions and basic results on independence algebras we refer the reader to \cite{ArEdGi,arfo,cameronSz,gould}. 

\begin{problem}
Let $A$ be a finite dimensional independence algebra. A pair $(a,G)$, where $a$ is a singular endomorphism of $A$ and  $G\leq\Aut(A)$, is said to be an $\Aut(A)$-pair if $$\langle a,G\rangle\setminus G=\langle a,\Aut(A)\rangle\setminus\Aut(A).$$

Classify  the $\mbox{Aut}(A)$-pairs.
\end{problem}	

This paper was prompted by the fact that the structure of the semigroups $\langle a,\sym\rangle\setminus \sym$ is well studied and known, as illustrated by Theorem \ref{illustr}. Therefore we propose the following problems.

\begin{problem}		
Let $G\leq \sym$ be a $2$-transitive group. (The list of those groups is available in \cite{cam,dixon}.)
For every $a\in \trans$ describe the structure of 
\[
\langle G,a\rangle\setminus G.
\]
In particular (where $G$ is a $2$-transitive group and $a\in \trans$):
\begin{enumerate}
\item  classify all the pairs $(a,G)$ such that $\langle a,G\rangle $ is a regular semigroup (that is, for all $x\in  \langle a,G\rangle$ there exists $y \in \langle a,G\rangle$ such that $x=xyx$); 
\item classify all the pairs $(a,G)$ such that $\langle a,G\rangle\setminus G $ is generated by its idempotents;
\item classify all the pairs $(a,G)$ such that $\langle a,G\rangle\setminus G =\langle g^{-1}ag\mid g\in G\rangle $; 
\item prove results analogous to those in Theorem \ref{illustr}.
\end{enumerate}
\end{problem}

\begin{problem}
For each $2$-transitive group $G$ classify the $G$-pairs, that is, the pairs $(a,H)$ such that  $H\leq \sym$, $a\in\trans$  and $\langle a,G\rangle\setminus G=\langle a,H\rangle\setminus H$.
\end{problem}

\begin{problem}	
Solve the analogues of the two previous problems for semigroups of linear transformations.
\end{problem}		

This paper deals with the general problem of showing how the group of units of a semigroup  shapes its structure. Regarding this general problem we recall here a number of very important and challenging problems (some of them related to automata theory). 

\begin{problem}	
\begin{enumerate}
\item Let $G$ be a Suzuki group $\Sz(q)$. Is it true that in the orbit of every  $3$-set contained in $\Omega$ there exists a section for every $3$-partition of $\Omega$?  The answer to this question  is necessary to decide if $G$ together with any rank $3$ map generates a regular semigroup (see \cite{ArCa12}). 
\item Classify the groups $G$ such that in the orbit of any $4$-subset there exists a section for every $4$-partition (both of $\Omega$), when   $\psl(2,q) \le G \le\pgaml(2,q)$, with either $q$ prime (except 
  $\psl(2,q)$ for $q\equiv1$ (mod~$4$)), or $q=2^p$ for $p$ prime. (See also \cite{ArCa12}). 
\item\label{synch} Classify the primitive groups $G\leq \sym$ such that there exists a partition $P$ of $N$ and a set $S\subseteq N$ such that $Sg$  is a section for $P$, for all $g\in G$ (see \cite{ben} and \cite{neu}).
\item Is there any primitive group $G\leq \sym$ and a non invertible transformation $a\in\trans\setminus \sym$ such that the kernel type of $a$ is non-uniform, but such that $\langle a,G\rangle$ does not generate a constant? (See \cite{ABC}.)
\item The diameter of a finite group is the maximum diameter of the group's Cayley graphs. Classify the primitive groups with diameter at most $n$ and that do not satisfy property (\ref{synch}). 
\end{enumerate}
\end{problem}


\begin{thebibliography}{9}

\bibitem{ABC}
J.~Ara\'ujo, Peter J. Cameron and W. Bentz.
\newblock{Groups Synchronizing a  Transformation of Non-Uniform Kernel.}
\newblock{http://arxiv.org/abs/1205.0682}

\bibitem{ArCa12}
J.~Ara\'ujo and Peter J. Cameron.
\newblock{Two Generalizations of Homogeneity in Groups with Applications to Regular Semigroups.}
\newblock{http://arxiv.org/pdf/1204.2195}

\bibitem{acmn}
J.~Ara\'ujo, Peter J. Cameron, James D. Mitchell, Max Neuhoffer.
\newblock{The classification of normalizing groups.}
\newblock{ {\em Journal of Algebra} \textbf{373} (2013), 481--490.}

\bibitem{ArEdGi}
J.~Ara\'ujo, M. Edmundo and S. Givant.
\newblock{$v^*$-Algebras, Independence Algebras and Logic.}
\newblock{{\em International Journal of Algebra and Computation} \textbf{21} (7) (2011), 1237--1257. }
 

\bibitem{arfo}
J.\ Ara\'{u}jo and J.\ Fountain.
 \newblock{The Origins of Independence
Algebras}
\newblock \textit{Proceedings of the Workshop on Semigroups and Languages
(Lisbon 2002)}, World Scientific, (2004), 54--67


\bibitem{arko}
J. Ara\'{u}jo and J. Konieczny,
\newblock{A method of finding automorphism groups of endomorphism monoids of relational systems.}
\newblock{{\em Discrete Math.} \textbf{307} 13, (2007), 1609--1620.} 

\bibitem{arko2}
J. Ara\'{u}jo and J. Konieczny,
\newblock{Automorphisms of endomorphism monoids of $1$-simple free algebras.}
\newblock{{\em  Comm. Algebra } \textbf{37} 1, (2009), 83--94.} 



\bibitem{ArMiSc}
J.~Ara\'ujo, J. D. Mitchell and C. Schneider.
\newblock{Groups that together with any transformation generate regular semigroup or idempotent generated semigroups.}
\newblock{ {\em Journal of Algebra} \textbf{343} (1) (2011), 93--106.}

\bibitem{ben}
F. Arnold and B. Steinberg.
\newblock{Synchronizing groups and automata.}
\newblock{ \textit{  Theoret. Comput. Sci.} \textbf{359} (2006), no. 1-3, 101--110.}

\bibitem{bp}
R.~A.~Beaumont and R.~P.~Peterson,
Set-transitive permutation groups,
\textit{Canad. J. Math.} \textbf{7} (1955), 35--42.

\bibitem{cam}
Peter~J. Cameron.
\newblock{{\em Permutation groups}, volume~45 of {\em London Mathematical
  Society Student Texts}.}
\newblock{Cambridge University Press, Cambridge, 1999.}

\bibitem{cameron-kantor}
Peter~J. Cameron and William~M. Kantor,
2-transitive and antiflag transitive collineation groups of finite projective spaces,
\textit{J. Algebra} \textbf{60} (1979), 384--422.


\bibitem{cameronSz}
P.\ J.\ Cameron and C.\ Szab\'{o},
\newblock{Independence algebras},
\newblock{{\em J.\ London Math.\ Soc.},  \textbf{61} (2000),  321--334.}

\bibitem{rjc}
Robin J. Chapman, Proof of Bertrand's Postulate, available from\hfil\break
\texttt{http://empslocal.ex.ac.uk/people/staff/rjchapma/etc/bertrand.pdf}

\bibitem{dixon}
John~D. Dixon and Brian Mortimer.
\newblock{{\em Permutation groups}, volume 163 of {\em Graduate Texts in
  Mathematics}.}
\newblock{Springer-Verlag, New York, 1996.}

\bibitem{mar}
O.  Ganyushkin and V.  Mazorchuk.
\newblock{Classical finite transformation semigroups. An introduction,}
 Algebra and Applications, 9. 
\newblock{Springer-Verlag London, Ltd., London, 2009.}
\bibitem{GAP}
The GAP~Group.
\newblock{{\em GAP -- Groups, Algorithms, and Programming, Version 4.6.2}},
  2013.
  \newblock{\verb+http://www.gap-system.org+}

\bibitem{gould}
V.\ Gould,  
\newblock{Independence algebras.}
 \newblock{{\em Algebra
Universalis}  \textbf{33} (1995), 294--318.}


\bibitem{hj}
R.~H\"aggkvist and J.~Janssen,
All-even latin squares,
\textit{Discrete Math.} \textbf{157} (1996), 199--206.

\bibitem{He74}
C. Hering.
\newblock{Transitive linear groups and linear groups which contain irreducible subgroups of prime order.}
\newblock{{\em Geometriae Dedicata} \textbf{2} (1974), 425--460.}

\bibitem{He85}
C. Hering.
\newblock{Transitive linear groups and linear groups which contain irreducible subgroups of prime order. II.}
\newblock{{\em J. Algebra} \textbf{93} (1985), 151--164.}

\bibitem{Ho95}
John~M. Howie.
\newblock{{\em Fundamentals of semigroup theory}, volume~12 of {\em London
  Mathematical Society Monographs. New Series}.}
\newblock{The Clarendon Press Oxford University Press, New York, 1995.}
\newblock{Oxford Science Publications.}

\bibitem{kantor}
W.~M.~Kantor, $4$-homogeneous groups,
\textit{Math. Z.} \textbf{103} (1968), 67--68.

\bibitem{kantor2}
W.~M.~Kantor, $k$-homogeneous groups,
\textit{Math. Z.} \textbf{124} (1972), 261--265.

\bibitem{kantor:inc}
William~M. Kantor.
\newblock{On incidence matrices of projective and affine spaces.}
\newblock{\textit{Math. Z.} \textbf{124} (1972), 315--318.}

\bibitem{kantor:line}
William~M. Kantor.
\newblock{Line-transitive collineation groups of finite projective spaces.}
\newblock{\textit{Israel J. Math.} \textbf{14} (1973), 229--235.}


\bibitem{levi96}
I. Levi.
\newblock{On the inner automorphisms of finite transformation semigroups.}
\newblock{\textit{Proc. Edinburgh Math. Soc. (2)}, \textbf{39} (1) (1996), 27--30.}

\bibitem{levi00}
I. Levi.
\newblock{Congruences on normal transformation semigroups.}
\newblock{\textit{Math. Japon. }, \textbf{52} (2) (2000), 247--261.}


\bibitem{lm}
I.~Levi and R.~B. McFadden.
\newblock{ {$S\sb n$}-normal semigroups.}
\newblock{{\em Proc. Edinburgh Math. Soc. (2)}, \textbf{37} (3) (1994), 471--476.}


\bibitem{lmm}
I.~Levi, D.~B. McAlister, and R.~B. McFadden.
\newblock{Groups associated with finite transformation semigroups.}
\newblock{{\em Semigroup Forum}, \textbf{61} (3) (2000), 453--467.}

\bibitem{Li86}
M.W. Liebeck.
\newblock{The affine permutation groups of rank 3.}
\newblock{{\em Bull. London Math. Soc.}, \textbf{18} (1986), 165--172.}


\bibitem{lw}
D.~Livingstone and A.~Wagner, Transitivity of finite permutation
groups on unordered sets,
\textit{Math. Z.} \textbf{90} (1965), 393--403.

\bibitem{lp}
T. {\L}uczak and L. Pyber,
On random generation of the symmetric group,
\textit{Combinatorics, Probability \& Computing} \textbf{2} (1993), 505--512.

\bibitem{ms}
W.~J.~Martin and B.~E.~Sagan, A new notion of transitivity for
groups and sets of permutations,
\textit{J. London Math. Soc} \textbf{73} (2006), 1--13.

\bibitem{mcalister}
Donald~B. McAlister.
\newblock{Semigroups generated by a group and an idempotent.}
\newblock{{\em Comm. Algebra}, \textbf{26} (2) (1998), 515--547.}

\bibitem{neu}
P. M. Neumann.
\newblock{Primitive permutation groups and their
section-regular partitions.}
\newblock{ \textit{Michigan Math. J.} {\bf 58} (2009), 309--322.}

\bibitem{vNM}
John von Neumann and Oskar Morgenstern,
\textit{Theory of Games and Economic Behavior},
Princeton University Press, Princeton, 1944.

\bibitem{Sc36}
J. Schreier, \"{U}ber Abbildungen einer abstrakten Menge Auf ihre
Teilmengen, {\it Fund.\ Math.\/} {\bf 28} (1936), 261--264.

\bibitem{symo}
J.S.V. Symons,
\textit{Normal transformation semigroups},
J. Austral. Math. Soc. Ser. A \textbf{22} (1976), no. 4, 385--390.

\bibitem{wie}
H. Wielandt, \textit{Finite Permutation Groups},
Academic Press, New York, 1964.

\end{thebibliography}
\end{document}